\documentclass[12pt]{amsart}
 \usepackage{mathtools}
\mathtoolsset{showonlyrefs,showmanualtags}

\usepackage[backrefs]{amsrefs}

 \usepackage[top=1.5in, bottom=1.5in, left=1.25in, right=1.25in]	{geometry}
\usepackage[all]{xy}
\usepackage{amsmath}
\usepackage{amssymb}
\usepackage{amsthm}
\usepackage{amscd}

\numberwithin{equation}{section}

\newtheorem{theorem}[equation]{Theorem}

\newtheorem{lemma}[equation]{Lemma}

\usepackage{hyperref} 
\hypersetup{
    colorlinks=true,       
    linkcolor=blue,          
    citecolor=magenta,        
    filecolor=magenta,      
    urlcolor=cyan           
}

\begin{document}
\title[Non-Linear Roth]{A Non-Linear Roth Theorem for Fractals of Sufficiently Large Dimension}
\author{Ben Krause}
\address{
Department of Mathematics,
Caltech \\
Pasadena, CA 91125}
\email{benkrause2323@gmail.com}
\date{\today}

\begin{abstract}
Suppose that $d \geq 2$, and that $A \subset [0,1]$ has sufficiently large dimension, $1 - \epsilon_d < \dim_H(A) < 1$. Then for any polynomial $P$ of degree $d$ with no constant term, there exists
$\{ x, x-t,x-P(t) \} \subset A$ with $t \approx_P 1$.
\end{abstract}

\maketitle


\section{Introduction}
In \cite{HLP}, the authors exhibit the existence of polynomial configurations in fractal sets; a key assumption on these fractal sets is that they have sufficiently large Fourier dimension, where we recall that the Fourier dimension of a set is given by
\[ \aligned
&\text{dim}_F(E) :=\\
& \qquad \sup\{ \beta : \text{E supports a probability measure}, \mu, \text{so that } | \hat{\mu}(\xi)| \leq C (1 + |\xi|)^{- \beta/2} \}. \endaligned\]

The purpose of this short note is to show that -- in one dimension -- this phenomenon is independent on Fourier dimension of fractal sets, provided that $E$ has sufficiently large Hausdorff dimension.

In particular, we have the following result.
\begin{theorem}\label{t:main}
Suppose that $A \subset [0,1]$ has sufficiently large dimension, $1 - \epsilon_d < \dim_H(A) < 1$, $d \geq 2$.
Then for any polynomial $P$ of degree $d$ with no constant term, there exists
$\{ x, x-t,x-P(t) \} \subset A$ with $t \approx_P 1$.
\end{theorem}

\subsection{Acknowledgement}
The author would like to thank Alex Iosevich for his perspective on the interplay between Fourier and Hausdorff dimension in detecting point configurations.

\subsection{Notation}
Here and throughout, $e(t) := e^{2\pi i t}$. For real numbers $A$ (typically taken to be dyadics), define $f_A$ to be the smooth Fourier restriction of $f$ to $|\xi| \approx A$, similarly define $f_{\leq A}$, etc.

For bump functions $\phi$, we let $\phi_j(x) := 2^j \phi(2^j x).$

We will make use of the modified Vinogradov notation. We use $X \lesssim Y$, or $Y \gtrsim X$, to denote the estimate $X \leq CY$ for an absolute constant $C$. We use $X \approx Y$ as shorthand for $Y \lesssim X \lesssim Y$. We also make use of big-O notation: we let $O(Y )$ denote a quantity that is $\lesssim Y$. If we need $C$ to depend on a parameter, we shall indicate this by subscripts, thus for instance $X \lesssim_p Y$ denotes the estimate $X \leq C_p Y$
for some $C_p$ depending on $p$. We analogously define $O_p(Y)$.

\section{The Argument}

\subsection{Preliminaries}
We need the following Lemmas; the first is essentially a consequence of the main result \cite{Broth}, or \cite{DGR}.

\begin{lemma}\label{l:roth}
Suppose that $0 \leq f_i = \phi_{t_i}*f$ for some bounded function $0 \leq f \leq \mathbf{1}_{[0,1]}$, \ $0 < t_i < \infty$ so that $\int f \geq \epsilon$. Then there exists $\delta_{\text{Roth}}(\epsilon) > 0$ so that
\[ \int \int f_1(x) f_2(x-t) f_3(x-P(t) ) \ dx dt \geq \delta_{\text{Roth}}(\epsilon).\]
\end{lemma}
The following refinement of \cite[Lemma 5]{Broth}, due to \cite{DGR}, is our primary tool.

\begin{lemma}[Lemma 1.4 of \cite{DGR}]\label{l:1.4}
Suppose that $\hat{f_2}$ is supported on $|\xi| \approx N$. Then there exists $\delta_0 > 0$ so that
\[ \left\| \int f_1(x-t) f_2(x-P(t)) \rho(t) \ dt \right\|_1 \lesssim N^{-\delta_0} \cdot \| f_1\|_2 \|f_2\|_2;\]
here $\rho = \rho_P$ is a $C^{\infty}$ function, supported in an annulus $\{|t| \approx_P 1\}$, with derivative $\| \partial^l \rho\|_{\infty} \lesssim_P 1$ for all $l \leq L$ sufficiently large.
\end{lemma}

Via the Fourier localization argument below, we see that this lemma essentially implies a Sobolev estimate. Since sets of dimension $>d_E$ support a probability measure $\mu$ with
\[ \| \phi_n * \mu \|_\infty \lesssim 2^{n(1-d_E - \kappa )},\]
for some $\kappa > 0$,
we see that $\mu$ is in the (negative) Sobolev class
\begin{equation}\label{e:Sob}
\mu \in H^{\kappa' - (1 - d_E)/2}
\end{equation}
for any $\kappa > \kappa' > 0$; this is essentially the key to the argument.

With these lemmas in hand, we turn to the proof.
\subsection{The Proof of Theorem \ref{t:main}}
\begin{proof}
Set $f = f_J := \mu* \phi_J$ for some sufficiently large $J$; it suffices to exhibit upper and lower bounds on
\begin{equation}\label{e:form}
    \int \int f(x) f(x-t) f(x-P(t)) \rho(t) \ dt dx
\end{equation}
independent of $J$, for $\rho$ an appropriate bump function.
Split \eqref{e:form} into two terms:
\begin{equation}\label{e:formLow}
    \int \int f(x) f(x-t) f_{\leq B} (x-P(t) ) \rho(t) \ dt dx
\end{equation}
and its complementary piece
\begin{equation}\label{e:formHigh}
    \int \int f(x) f(x-t) f_{> B} (x-P(t)) \rho(t) \ dt dx
\end{equation}
where $B$ is a large parameter to be determined later.\footnote{In particular, we will choose $B$ so large that for $c \geq c_0 > 0$ bounded away from zero, $\delta_{\text{Roth}}(c) \gg B^{-\delta_0}$, where $\delta_0$ is as in Lemma \ref{l:1.4}.}

We begin with \eqref{e:formLow}, which we write as
\[ \int \int \hat{f}(-\xi - \eta) \hat{f}(\xi) \widehat{f_{\leq B}}(\eta) m(\xi,\eta) \ d\xi d\eta,\]
where
\[ \aligned
m(\xi,\eta) &:= \int e(-\xi t - \eta P(t)) \rho(t) \ dt \\
& = m(\xi,\eta) \cdot \mathbf{1}_{|\xi| \lesssim_P |\eta|} + \mathbf{1}_{|\xi| \gg_P |\eta|} \cdot O_B((1+|\xi|)^{-B}) \\
& \qquad =: m_1(\xi,\eta) + m_2(\xi,\eta)
\endaligned\]
by non-stationary phase considerations.

Now, with $B \gtrsim A = A(B,P) \gg B$ a large threshold, decompose \eqref{e:formLow} as a sum of three terms:

\begin{align}
\eqref{e:formLow} &=
\int \int f_{\lesssim A} (x) f_{\leq A} (x-t) f_{\leq B} (x-P(t)) \rho(t) \ dt dx
\label{e:1}\\
& \qquad + \int \int f_{\gg A} (x) f_{\leq A} (x-t) f_{\leq B} (x-P(t)) \rho(t) \ dt dx
\label{e:2}\\
& \qquad \qquad \int \int f(x) f_{> A} (x-t) f_{\leq B} (x-P(t)) \rho(t) \ dt dx.
\label{e:3}
\end{align}
The first term is a main term;
an upper bound is simply given by
\[ \eqref{e:1} \lesssim \| f_{\lesssim A} \|_\infty \cdot
\| f_{\leq A} \|_\infty \cdot
\| f_{\leq B} \|_\infty \lesssim
A^{3(1- d_E)};\]
as we will see, \eqref{e:2} and \eqref{e:3} are lower order error terms, so this upper bound majorizes \eqref{e:formLow}.

As for the lower bound, an application of Lemma \ref{l:roth} yields a lower bound of
\[ \eqref{e:1} \gtrsim \delta_{\text{Roth}}(A^{d_E-1});\]
the loss of $A^{d_E -1}$ comes from reproducing: we have
\[ \|f_{\lesssim A}\|_\infty \lesssim A^{1-d_E}.\]
Since $\delta_{\text{Roth}}(\epsilon)$ grows super-polynomially in $\epsilon^{-1}$, see \cite{SA}, we stipulate that 
\begin{equation}\label{e:d}
1 - \frac{C}{\log B} < d_E < 1.
\end{equation}

The second term, \eqref{e:2}, vanishes identically, since
\[ \widehat{f_{\gg A}}(-\xi - \eta) \widehat{f_{\leq A}}(\xi) \widehat{f_{\leq B}}(\eta) = 0.\]
We express the third term using the Fourier transform:
\[ \aligned \eqref{e:3} &=
\int \int \hat{f}(-\xi - \eta) \widehat{f_{>A}}(\xi) \widehat{f_{\leq B}}(\eta) m_1(\xi,\eta) \ d\xi d\eta \\
& \qquad + \int \int \hat{f}(-\xi - \eta) \widehat{f_{>A}}(\xi) \widehat{f_{\leq B}}(\eta) m_2(\xi,\eta)  \ d\xi d\eta.
\endaligned \]
The first term vanishes identically since $|\xi|$ is so much larger that $|\eta|$, for an appropriate choice of $A$.
As for the error term, we estimate:
\[ \aligned
& \left| \int \int \hat{f}(-\xi - \eta) \widehat{f_{>A}}(\xi) \widehat{f_{\leq B}}(\eta) m_2(\xi,\eta) \ d\xi d\eta \right| \\
& \qquad \leq \int \int |\widehat{f_{>A}}(\xi)| \cdot |\widehat{f_{\leq B}}(\eta)| \cdot |m_2(\xi,\eta)| \ d\xi d\eta \\
& \qquad \qquad \leq \sum_{N > A} \int \int  |\widehat{f_{N}}(\xi)| \cdot  |\widehat{f_{\leq B}}(\eta)| \cdot |\xi|^{-B} \ d\xi d\eta \\
& \qquad \qquad \qquad \lesssim B^{1 - d_E/2} \cdot A^{1-d_E/2 - B}  \lesssim A^{2-d_E - B} =: A^{-C};
\endaligned\]
note the use of the trivial estimate $\| \hat{f} \|_\infty \leq \|f\|_1$ in passing to the second line.

In particular, taking into account \eqref{e:d}, we have exhibited
\begin{equation}\label{e:formLfin}
    \eqref{e:formLow} \gtrsim \delta_{\text{Roth}}(A^{d_E -1 }) - A^{-C} \gtrsim \delta_{\text{Roth}}(c)
\end{equation}
for some $c \geq c_0 > 0$ bounded away from zero, provided that we have chosen $A$ sufficiently large.

We next term to \eqref{e:formHigh}, which we decompose as a sum of $N > B$:

\begin{align}
\eqref{e:formLow} &=
\int \int f (x) f (x-t) f_{> B} (x-P(t) ) \rho(t) \ dt dx
\label{ee:1}\\
& \qquad = \sum_{N > B} \int \int \widehat{f_{\lesssim N}}(-\xi-\eta) \widehat{f_{\lesssim N}}(\xi) \widehat{f_N}(\eta) m_1(\xi,\eta) \ d\xi d\eta \label{ee:2} \\
& \qquad \qquad
+
\sum_{N > B} \int \int \widehat{f_{\gg N}}(-\xi-\eta) \widehat{f_{\gg N}}(\xi) \widehat{f_N}(\eta) m_2(\xi,\eta) \ d\xi d\eta \label{ee:3}
\end{align}
We begin with \eqref{ee:3}; by arguing as previously, the $N$th term admits an upper bound of $N^{-C}$ for a very large $C = C(B)$, which leads to the estimate
\[ |\eqref{ee:3}| \lesssim B^{-C}. \]
It remains to consider \eqref{ee:2}; we extract the $N$th term once again,
\[ \aligned
& \left| \int \int \widehat{f_{\lesssim N}}(-\xi-\eta) \widehat{f_{\lesssim N}}(\xi) \widehat{f_N}(\eta) m_1(\xi,\eta) \ d\xi d\eta \right| \\
& \qquad = \left| \int \int \widehat{f_{\lesssim N}}(-\xi-\eta) \widehat{f_{\lesssim N}}(\xi) \widehat{f_N}(\eta) m (\xi,\eta) \ d\xi d\eta \right| \\
& \qquad \qquad \leq \|f_{\lesssim N}\|_\infty \cdot \| \int f_{\lesssim N}(x-t) f_{N}(x- P(t) ) \rho(t) \ dt \|_{L^1}  \\
& \qquad \qquad \qquad \lesssim N^{1 - d_E} \cdot N^{-\delta_0} \cdot \| f_{\lesssim N} \|_2^2 \\
& \qquad \qquad \qquad \qquad \lesssim N^{2 - 2d_E - \delta_0}. \endaligned\]
In passing from the first line to the second, we have (possibly) discarded $O(1)$ terms of the form
\[ \int \int |\widehat{f_{\approx C N}}(-\xi-\eta)| \cdot  | \widehat{f_{\approx C' N}}(\xi)| \cdot | \widehat{f_N}(\eta) | \cdot |m_1(\xi,\eta)|  \ d\xi d\eta \]
for some large $1 \ll C, C' \lesssim 1$ as we drop the Fourier restriction in the definition of $m_1$; but, on this domain, we retain the pointwise bound $|m_1(\xi,\eta)| \lesssim (1 + |\xi|)^{-B}$, so we may handle these error terms as above.

In particular, since \eqref{e:d} ensures that we have $1 - \delta_0/2 < d_E$ for sufficiently large $B$, we have exhibited a upper bound
\begin{equation}\label{e:formHfin}
    \eqref{e:formHigh} \lesssim B^{- \delta_0 - 2d_E + 2} \lesssim B^{-\delta_0}
\end{equation}

Combining \eqref{e:formLfin} and \eqref{e:formHfin}, we see that we may estimate from below
\[ \eqref{e:form} \geq \delta_{\text{Roth}}(c) - C B^{-\delta_0},\]
which yields the result.
\end{proof}

\typeout{get arXiv to do 4 passes: Label(s) may have changed. Rerun}

\end{document}